\documentclass[12pt]{article}
\usepackage{latexsym,amssymb}
\usepackage{amsmath}
\usepackage[mathscr]{eucal}
\usepackage{color}
\usepackage{enumerate}

\usepackage{theorem}
\newtheorem{theorem}{Theorem}[section]
\newtheorem{proposition}[theorem]{Proposition}

\newtheorem{corollary}[theorem]{Corollary}

\theorembodyfont{\rmfamily}

\newcommand{\supp}{{\rm supp }}
\newcommand{\vol}{{\rm vol }}

\newcommand{\R}{{\mathbb{R}}}
\numberwithin{equation}{section}

\newcommand{\C}{\mathcal{C}}
\newcommand{\cP}{\mathcal{P}}
\newcommand{\CD}{\mathrm{CD}}
\newcommand{\RCD}{\mathrm{RCD}}
\newcommand{\Ric}{\mathrm{Ric}}
\newcommand{\Hess}{\mathrm{Hess}}
\newcommand{\e}{\mathrm{e}}
\newcommand{\ac}{\mathrm{ac}}

\begin{document}

\title{Equality in the logarithmic Sobolev inequality}

\author{Shin-ichi OHTA\thanks{
Department of Mathematics, Osaka University, Osaka 560-0043, Japan
({\sf s.ohta@math.sci.osaka-u.ac.jp})} \textsuperscript{,}\footnotemark[3]
\and
Asuka TAKATSU\thanks{
Department of Mathematical Sciences, Tokyo Metropolitan University, Tokyo 192-0397, Japan
({\sf asuka@tmu.ac.jp})} \textsuperscript{,}\thanks{
RIKEN Center for Advanced Intelligence Project (AIP),
1-4-1 Nihonbashi, Tokyo 103-0027, Japan}}

\date{\today}

\maketitle

\vspace{-5mm}
\begin{center}
{\it  Dedicated to the memory of Kazumasa Kuwada}
\end{center}
\smallskip

\begin{abstract}
We investigate the rigidity problem for the logarithmic Sobolev inequality
on weighted Riemannian manifolds satisfying $\Ric_{\infty} \ge K>0$.
Assuming that equality holds,
we show that the $1$-dimensional Gaussian space is necessarily split off,
similarly to the rigidity results of Cheng--Zhou on the spectral gap
as well as Morgan on the isoperimetric inequality.
The key ingredient of the proof is the needle decomposition method
introduced on Riemannian manifolds by Klartag.
We also present several related open problems.
\end{abstract}

\section{Introduction}
The rigidity and stability properties of various inequalities
are one of the most fundamental problems in geometric analysis,
and recently draw considerable attention
from the geometric, analytic and probabilistic viewpoints.
We refer to \cite{CFMP,MN,El,BBJ} (resp.\ \cite{CMM})
for the quantitative isoperimetry on Gaussian spaces (resp.\ positively curved spaces),
and to \cite{CZ,Ma1,GKKO} for rigidity results on the sharp spectral gap
(equivalently, Poincar\'e inequality), to name a few.

In this article, we study the rigidity problem for the logarithmic Sobolev inequality
on positively curved weighted Riemannian manifolds.
On a Riemannian manifold $(M,g)$
equipped with a weighted probability measure $\omega=\e^{-\Psi} \,\vol_M$
($\Psi \in \C^{\infty}(M)$)
and satisfying the weighted Ricci curvature bound
$\Ric_{\infty}:=\Ric_g +\Hess_g \Psi \ge K>0$,
the \emph{logarithmic Sobolev inequality}
\begin{equation}\label{eq:LSI}
\int_M \rho \log \rho \,d\omega
 \leq \frac{1}{2K} \int_M \frac{|\nabla \rho|^2}{\rho} \,d\omega
\end{equation}
holds for all nonnegative locally Lipschitz functions $\rho: M \longrightarrow [0,\infty)$
with $\|\rho\|_{L^1(\omega)}=1$ (see Section~\ref{sc:prel} for a further explanation).

We will characterize when equality holds in \eqref{eq:LSI}.
We remark that, since the logarithmic Sobolev inequality implies the Poincar\'e inequality,
equality in \eqref{eq:LSI} is a weaker assumption than that for the spectral gap.
In the Gaussian case ($M=\R^n$, $d\omega=(K/2\pi)^{n/2}\e^{-K|x|^2/2} \,dx$),
Carlen \cite[Theorem~4]{Ca} showed that equality holds if and only if $\rho \omega$
is a translation of $\omega$ (see also \cite{Le}).
This was generalized to weighted Euclidean spaces satisfying $\Ric_{\infty} \ge K$
in \cite[Theorem~3.11]{AMTU} (see also \cite[p.~390]{OV}).
The manifold case was briefly discussed in \cite[p.~391]{OV},
and our main theorem below gives a precise answer to this question.

\begin{theorem}[Main theorem]\label{main}
Let $(M,g,\omega)$ be a connected, complete,
weighted $\C^{\infty}$-Riemannian manifold of dimension $n \ge 2$
such that $\Ric_{\infty} \ge K>0$.
Assume that there is a nonconstant, nonnegative, locally Lipschitz function
$\rho: M \longrightarrow [0,\infty)$ such that $\|\rho\|_{L^1(\omega)}=1$ and 
\[ \int_M \rho \log \rho\, d\omega =
 \frac{1}{2K} \int_M \frac{|\nabla \rho|^2}{\rho}\, d\omega. \]
Then we have the following.

\begin{enumerate}[{\rm (i)}]
\item\label{main1}
There exists a weighted Riemannian manifold $(N,g_N,\omega_N)$ of dimension $n-1$
satisfying $\Ric_{\infty} \ge K$ again such that
$(M,g,\omega)$ is isometric to the product of
$(\mathbb{R},|\cdot|,\gamma_K)$ and $(N,g_N,\omega_N)$ as weighted Riemannian manifolds,
where $|\cdot|$ denotes the Euclidean metric and
$\gamma_K$ denotes the Gaussian measure with variance $K^{-1}$, namely 
\[ d\gamma_K :=\sqrt{\frac{K}{2\pi}} \e^{-Kx^2/2} \,dx. \]

\item\label{main2}
In the product structure as in \eqref{main1},
$\rho$ is a translation of the density of $\omega$ in the $\R$-direction.
Precisely, there exists some $t \in \R \setminus \{0\}$ such that
the isometry $T_t:\R \times N \longrightarrow \R \times N$ defined by $T_t(s,x):=(s+t,x)$
pushes $\omega$ forward to $\rho \omega$.
\end{enumerate}
\end{theorem}

The first assertion (\ref{main1}) has the same consequence
as Cheng--Zhou's result on the sharp spectral gap \cite[Theorem~2]{CZ}
as well as Morgan's result on the sharp isoperimetric inequality \cite[Theorem~18.7]{Mo}
under $\Ric_{\infty} \ge K$.
We will in fact show that equality in \eqref{eq:LSI} implies
the sharp spectral gap and then apply Cheng--Zhou's result
(in the same spirit as the argument in \cite{Ma2} concerning equality of isoperimetric inequalities).
The second assertion (\ref{main2}) shows a similar phenomenon to the characterizations in \cite{Ca,AMTU}.
The key tool in the proof of Theorem~\ref{main} is the \emph{needle decomposition}
established by Klartag \cite{Kl} (see Theorem~\ref{needle}),
we refer to \cite{CM1,CM2,CMM,Ma2} for other applications of this powerful method.

Combining Theorem~\ref{main} with the results in \cite{Mo,CZ},
we complete the following chain of characterizations of the sharp geometric and functional inequalities.

\begin{corollary}\label{cr:main}
Let $(M,g,\omega)$ be a connected, complete, weighted $\C^{\infty}$-Riemannian manifold
of dimension $n \ge 2$ such that $\Ric_{\infty} \ge K>0$.
Then the following are equivalent.

\begin{enumerate}[{\rm (I)}]
\item
Equality holds in the logarithmic Sobolev inequality \eqref{eq:LSI} on $(M,g,\omega)$.

\item
$(M,g,\omega)$ attains the sharp spectral gap $\lambda_1=K$ or, equivalently,
equality holds in the Poincar\'e inequality \eqref{eq:Poin}.

\item
Equality holds in the Bakry--Ledoux isoperimetric inequality in $(M,g,\omega)$.

\item
$(M,g,\omega)$ splits off the $1$-dimensional Gaussian space of variance $K^{-1}$
in the sense of Theorem~$\ref{main}$\eqref{main1}.
\end{enumerate}
\end{corollary}

The article is organized as follows.
In the next section, we recall key definitions and tools necessary in our argument,
including the curvature-dimension condition and needle decompositions,
from optimal transport theory.
Section~\ref{sc:proof} is devoted to the proof of Theorem~\ref{main}.
We finally discuss several further problems in Section~\ref{outro}.

\section{Preliminaries for optimal transport theory}\label{sc:prel}

Throughout the article,
let $(M,g)$ be a connected, complete $\C^{\infty}$-Riemannian manifold without boundary
of dimension $n \geq 2$.
We will denote by $d_M$ and $\vol_M$ the Riemannian distance function
and the volume measure on $M$, respectively.
Let us also fix a probability measure $\omega:=\e^{-\Psi} \,\vol_M$ on $M$ with $\Psi \in \C^\infty(M)$.
Then $(M,g,\omega)$ is called a \emph{weighted Riemannian manifold}
or a \emph{Riemannian manifold with density}.

Optimal transport theory is quite useful in geometry and analysis of weighted manifolds.
In this section, we recall some necessary notions and results.
We refer to \cite{Vi} for the fundamentals of optimal transport theory
as well as the \emph{curvature-dimension condition},
and to \cite{Kl,CM1} for the recent breakthrough
called the \emph{needle decomposition} or the \emph{localization}.

\subsection{Wasserstein distance}

Let $\cP_2(M)$ be the set of all Borel probability measures $\mu$ on $M$
with finite second moments, namely $d_M(o,\cdot) \in L^2(\mu)$
for some (and hence any) $o\in M$.
We denote by $\cP_2^{\ac}(M,\omega)$ the subset of $\cP_2(M)$
consisting of absolutely continuous measures with respect to $\omega$.

Given $\mu_0,\mu_1 \in \cP_2(M)$,
we call $\pi \in \cP(M \times M)$ a \emph{coupling} of $(\mu_0,\mu_1)$ if
\[
\pi(A \times M)=\mu_0(A), \qquad \pi(M \times A)=\mu_1(A)
\]
hold for any Borel set $A \subset M$.
In other words, the push-forward measure of $\pi$ to the first (resp.\ second) component $M$
coincides with $\mu_0$ (resp.\ $\mu_1$).
We denote by $\Pi(\mu_0, \mu_1)$ the set of all couplings of $(\mu_0,\mu_1)$.
Then the \emph{$L^2$-Wasserstein distance} between $\mu_0,\mu_1 \in \cP_2(M)$ is defined as 
\[
W_2(\mu_0,\mu_1):=\inf_{\pi \in \Pi(\mu_0,\mu_1)}  \| d_M\|_{L^2(\pi)}.
\]
A coupling $\pi \in \Pi(\mu_0, \mu_1)$ is said to be \emph{optimal} if it attains the infimum above.
For any $\mu_0,\mu_1 \in \cP_2(M)$, $W_2(\mu_0,\mu_1)$ is finite and  
an optimal coupling always exists (see \cite[Chapter~4]{Vi}).
Moreover, the pair $(\cP_2(M),W_2)$ is a complete separable metric space
(see \cite[Chapter~6]{Vi}).

In the case where $\mu_0 \in \cP_2^{\ac}(M,\omega)$,
an optimal coupling is given by the push-forward measure of $\mu_0$ as follows
(see \cite{Br,Mc,FG,Vi}).

\begin{theorem}[Optimal transports]\label{th:FG}
Given any $\mu_0 \in \cP_2^{\ac}(M,\omega)$ and $\mu_1 \in \cP_2(M)$,
there exists a locally semi-convex function $\phi:\Omega \to \mathbb{R}$
on an open set $\Omega \subset M$ with $\mu_0(\Omega)=1$
such that the map  $T_t(x):=\exp_x(t\nabla\phi(x))$, $t \in [0,1]$,
provides a unique minimal geodesic from $\mu_0$ to $\mu_1$ with respect to $W_2$.
Precisely, $(T_0 \times T_1)_{\sharp}\mu_0$ is a unique optimal coupling of $(\mu_0,\mu_1)$,
and $\mu_t :=(T_t)_{\sharp}\mu_0$ gives the unique $W_2$-minimal geodesic from $\mu_0$ to $\mu_1$.
Moreover, $\mu_t \in \cP_2^{\ac}(M,\omega)$ holds for all $t \in [0,1)$.
\end{theorem}

Denoted by $(T_t)_{\sharp}\mu_0$ is the push-forward measure of $\mu_0$ by $T_t$.
A locally semi-convex function is locally Lipschitz and twice differentiable almost everywhere
by the Alexandrov--Bangert theorem (see \cite{Ba,Vi}).
Hence $T_t$ is differentiable almost everywhere on $\Omega$,
and $\mu_t=(T_t)_{\sharp}\mu_0$ implies the following \emph{Monge--Amp\`ere equation}
(or the \emph{Jacobian equation})
\begin{equation}\label{MAonM}
\rho_t\big( T_t(x) \big) \cdot \e^{\Psi(x)-\Psi(T_t(x))} {\det}_g \big( dT_t(x) \big) =\rho_0(x)
\end{equation}
for all $t \in [0,1)$ at $\mu_0$-almost every $x \in \Omega$,
where $\mu_t=\rho_t \omega$ and $\det_g(dT_t)$ is the Jacobian determinant of $T_t$ with respect to $g$.
If $\mu_1 \in \cP_2^{\ac}(M,\omega)$, then \eqref{MAonM} holds also at $t=1$.

\subsection{Curvature-dimension condition}

The deep relation between optimal transport theory and Riemannian geometry
was revealed in the \emph{curvature-dimension condition},
the convexity of the relative entropy characterizing the lower weighted Ricci curvature bounds.

For $\mu=\rho \omega \in \cP_2^{\ac}(M,\omega)$,
we define the \emph{relative entropy} with respect to $\omega$ as
\[
H_\omega(\mu) :=\int_M \rho \log \rho\, d\omega.
\]
We need to modify the Ricci curvature of $(M,g)$,
taking the effect of the weight function $\Psi$ into account,
to the \emph{weighted Ricci curvature}
(also called the \emph{Bakry--\'Emery--Ricci curvature})
\[ \Ric_{\infty}(v) :=\Ric_g(v,v) +\Hess_g\Psi(v,v) \]
for $v \in TM$,
where $\Ric_g$ is the Ricci curvature tensor and $\Hess_g \Psi$ is the Hessian of $\Psi$
with respect to $g$.
We will say that $\Ric_{\infty} \ge K$ holds for some $K \in \R$
if we have $\Ric_{\infty}(v) \ge K|v|^2$ for all $v \in TM$.
The next theorem was established in \cite{CMS1,CMS2,vRS,Stcon}.

\begin{theorem}[Curvature-dimension condition]\label{th:CD}
For a weighted Riemannian manifold $(M,g,\omega)$ and $K\in \mathbb{R}$,
we have $\Ric_{\infty} \ge K$ if and only if $H_{\omega}$ is \emph{$K$-convex} in the sense that,
for any pair $\mu_0, \mu_1 \in \cP_2^{\ac}(M,\omega)$,
\begin{equation}\label{eq:CD}
H_{\omega}(\mu_t) \leq (1-t) H_\omega(\mu_0) +t H_\omega(\mu_1)
 -\frac{K}{2}t(1-t)W_2(\mu_0,\mu_1)^2
\end{equation}
holds for all $t \in (0,1)$, where $(\mu_t)_{t \in [0,1]} \subset \cP_2^{\ac}(M,\omega)$
is the unique minimal geodesic from $\mu_0$ to $\mu_1$ as in Theorem~$\ref{th:FG}$.
\end{theorem}

The above $K$-convexity of $H_{\omega}$ with respect to $W_2$ is called
the \emph{curvature-dimension condition} $\CD(K,\infty)$.
We refer to \cite{StI,StII,LV,Vi} and recent \cite{Am} for the rapidly developing theory
of metric measure spaces satisfying the curvature-dimension condition.

Going back to the pioneering work \cite{OV},
it is known that the curvature-dimension condition intimately connects
the curvature bound $\Ric_{\infty} \ge K$ with some functional inequalities.
For $\mu=\rho \omega \in \cP_2^{\ac}(M,\omega)$ with locally Lipschitz density $\rho$,
the \emph{Fisher information} with respect to $\omega$ is defined as
\[ I_\omega(\mu) :=\int_M \frac{|\nabla \rho|^2}{\rho}\, d\omega, \]
where we set $|\nabla \rho|^2/\rho :=0$ on $\rho^{-1}(0)$.
The Fisher information appears as the derivative of $H_{\omega}$ along $W_2$-geodesics,
and then the $K$-convexity \eqref{eq:CD} shows the following
(see \cite{OV}, \cite[\S 6]{LV}).

\begin{proposition}[Talagrand, HWI inequalities]\label{pr:HWI}
If $(M,g,\omega)$ satisfies $\Ric_{\infty} \ge K>0$, then we have
\begin{equation}\label{Talagrand}
W_2(\mu,\omega)^2 \le \frac{2}{K} H_\omega(\mu)
\end{equation}
for any $\mu \in \cP_2^{\ac}(M,\omega)$, and
\begin{equation}\label{eq:HWI}
H_\omega(\mu) \le W_2(\mu,\omega) \sqrt{I_\omega(\mu)} -\frac{K}{2} W_2(\mu,\omega)^2
\end{equation}
for any $\mu \in \cP_2^{\ac}(M,\omega)$ whose density function is locally Lipschitz.
\end{proposition}

The first inequality \eqref{Talagrand} is called the \emph{Talagrand inequality}
and \eqref{eq:HWI} is called the \emph{HWI inequality}.
Completing the square in \eqref{eq:HWI},
we obtain the logarithmic Sobolev inequality \eqref{eq:LSI} as
\begin{equation}\label{HWItoLSI}
H_\omega(\mu)
 \le \frac{1}{2K} I_\omega(\mu)
 -\frac{K}{2}\bigg( \frac{1}{K} \sqrt{I_{\omega}(\mu)}-W_2(\mu,\omega) \bigg)^2
 \le \frac{1}{2K} I_\omega(\mu).
\end{equation}
Moreover, we can deduce from the logarithmic Sobolev inequality the \emph{Poincar\'e inequality}
\begin{equation}\label{eq:Poin}
\int_M f^2 \,d\omega -\bigg( \int_M f \,d\omega \bigg)^2
 \le \frac{1}{K} \int_M |\nabla f|^2 \,d\omega
\end{equation}
for any locally Lipschitz function $f \in L^2(\omega)$ (see \cite{OV}, \cite[\S 6]{LV}).

\subsection{Needle decompositions}

Let us close the section by recalling the \emph{needle decomposition}
(also called the \emph{localization}).
This is a powerful method, with its root in convex geometry,
for reducing a high-dimensional problem to its $1$-dimensional counterpart.
We refer to \cite{Kl} for the background and the generalization to weighted Riemannian manifolds,
and to \cite{CM1} for a further generalization to metric measure spaces satisfying
the curvature-dimension condition.
We will use the following version in \cite[Theorems~1.2, 1.5]{Kl}.

\begin{theorem}[Needle decompositions]\label{needle}
Assume that $(M,g,\omega)$ satisfies $\Ric_{\infty} \ge K$
and let $f \in L^1(\omega)$ satisfy
\[ \int_M f \,d\omega=0 \]
and $f \cdot d_M(o,\cdot) \in L^1(\omega)$ for some $($and hence all$)$ $o \in M$.
Then there exists a $1$-Lipschitz function $u$,
a partition $L$ of $M$, a probability measure $\nu$ on $L$
and a family of probability measures $\{\omega^{\ell}\}_{\ell\in L}$ on $M$
satisfying the following conditions.
\begin{enumerate}[{\rm (i)}]
\item\label{needle1}
For any measurable set $A \subset M$, we have
\[ \omega(A)=\int_L \omega^{\ell}(A) \,\nu(d\ell). \]
\item\label{needle2}
For $\nu$-almost all $\ell \in L$, we have $\supp (\omega^{\ell}) \subset \ell$
and $\ell$ is a \emph{transport ray} associated with $u$.
Moreover, if $\ell$ is not a singleton,
then $(\ell, |\cdot|, \omega^{\ell})$ satisfies $\Ric_{\infty} \ge K$ in the weak sense.
\item\label{needle3}
For $\nu$-almost all $\ell\in L$, 
we have
\[ \int_\ell f \,d\omega^{\ell}=0. \]
\end{enumerate}
\end{theorem}

The $1$-Lipschitz function $u$ is called the \emph{guiding function}.
In (\ref{needle2}), $\ell \subset M$ being a \emph{transport ray} means that
$|u(x)-u(y)|=d_M(x,y)$ holds for any $x,y \in \ell$
and $\ell$ cannot be extended keeping this property.
A transport ray is necessarily the image of a minimal geodesic,
and $|\cdot|$ denotes the distance structure of $d_M$ restricted on it.
Also in (\ref{needle2}), $\Ric_{\infty} \ge K$ in the weak sense means that,
if we write $d\omega^{\ell}=\e^{-\Psi_{\ell}} \,dx$ (by identifying $\ell$ with an interval in $\R$),
then $\Psi_{\ell}$ is $K$-convex in the same way as \eqref{eq:CD}.
The construction of the needle decomposition is closely related to the $L^1$-optimal transport theory.
In such a context, the function $f$ is regarded as the difference between two probability measures,
the guiding function $u$ is a Kantorovich potential,
and the $L^1$-optimal transport is done along \emph{needles} $\ell$.

\section{Proof of Theorem~\ref{main}}\label{sc:proof}

Let $(M,g,\omega)$ be a weighted Riemannian manifold as in Theorem~\ref{main},
namely it satisfies $\Ric_{\infty} \ge K>0$ and
there is a nonconstant, nonnegative, locally Lipschitz function $\rho$ with $\|\rho\|_{L^1(\omega)}=1$
attaining equality in the logarithmic Sobolev inequality \eqref{eq:LSI}.
Notice that $\mu:=\rho\omega \in \cP_2(M)$ by the Talagrand inequality \eqref{Talagrand}
(precisely, we take the truncation $\rho_i:=\min\{\rho,i\}$ and apply \eqref{Talagrand}
to the normalization of $\rho_i$, then we have \eqref{Talagrand} for $\mu$ as $i \to \infty$
thanks to the lower semi-continuity of $W_2$ found in \cite[Lemma~4.3]{Vi} for instance).
Recalling \eqref{HWItoLSI}, then we have
\begin{equation}\label{HWIeq}
H_\omega(\mu) = \frac{1}{2K} I_\omega(\mu)=\frac{K}{2}W_2(\mu,\omega)^2.
\end{equation}
We shall show that $\log\rho$ attains equality in the Poincar\'e inequality \eqref{eq:Poin}.
Then the rigidity result in \cite{CZ} applies.

\subsection{Decomposition into needles}\label{ssc:decomp}

We first recall  how to derive the logarithmic Sobolev inequality on $M$
from those on needles (see \cite{CM2}).
Put $f:=\rho-1$ and notice that $f \in L^1(\omega)$ and $\int_M f \,d\omega=0$.
Moreover, since $\mu \in \cP_2(M)$, we have for any $o \in M$
\begin{align*}
&\int_M |f(x)| d_M(o,x) \,\omega(dx)
 \le \int_M \big( \rho(x)+1 \big) d_M(o,x) \,\omega(dx) \\
&\le \bigg( \int_M \rho(x) \,\omega(dx) \bigg)^{1/2}
 \bigg( \int_M \rho(x) d_M(o,x)^2 \,\omega(dx) \bigg)^{1/2}
 +\int_M d_M(o,x) \,\omega(dx) \\
&<\infty.
\end{align*}
In the last inequality, we used the fact that $\omega$ has the Gaussian decay
(\cite[Theorem~4.26]{StI}) to see $\int_M d_M(o,x) \,\omega(dx)<\infty$.
Hence we can apply Theorem~\ref{needle} to $f$
and obtain a $1$-Lipschitz function $u$,
a partition $L$ of $M$ along with a family of probability measures $\{\omega^{\ell}\}_{\ell \in L}$,
and a probability measure $\nu$ on $L$.

Theorem~\ref{needle}(\ref{needle3}) shows $\int_\ell f \,d\omega^{\ell}=0$ $\nu$-almost every $\ell \in L$,
which yields $\int_\ell \rho \,d\omega^{\ell}=1$.
Thus, since $(\ell,|\cdot|,\omega^{\ell})$ enjoys $\Ric_{\infty} \ge K$ by Theorem~\ref{needle}(\ref{needle2}),
the $1$-dimensional logarithmic Sobolev inequality
on the needle $(\ell,|\cdot|,\omega^{\ell})$ yields
\begin{equation}\label{1dimLSI}
\int_{\ell} \rho \log \rho \,d\omega^{\ell}
 \le \frac{1}{2K} \int_{\ell} \frac{|\nabla^{\ell} \rho|^2}{\rho} \,d\omega^{\ell},
\end{equation}
where $|\nabla^{\ell} \rho|$ denotes the slope of $\rho$ taken along $\ell$.
Note that $|\nabla^{\ell} \rho| \le |\nabla \rho|$,
thereby the integration of \eqref{1dimLSI} in $\nu$ provides the logarithmic Sobolev inequality on $M$
(by Theorem~\ref{needle}(\ref{needle1})).
Under our assumption \eqref{HWIeq}, we have equality in \eqref{1dimLSI} and
\begin{equation}\label{1dim-eq}
\int_{\ell} \rho \log \rho \,d\omega^{\ell}
 =\frac{1}{2K} \int_{\ell} \frac{|\nabla^{\ell} \rho|^2}{\rho} \,d\omega^{\ell}
 =\frac{1}{2K} \int_{\ell} \frac{|\nabla \rho|^2}{\rho} \,d\omega^{\ell}
 =\frac{K}{2}W_2(\mu^{\ell},\omega^{\ell})^2
\end{equation}
$\nu$-almost every $\ell \in L$, where $\mu^{\ell}:=\rho \omega^{\ell} \in \cP_2^{\ac}(\ell,\omega^{\ell})$.

\subsection{Analysis on needles}\label{ssc:1dim}

Next, in order to consider the outcome of the equality \eqref{1dim-eq},
let us fix a needle $(\ell,|\cdot|,\omega^{\ell})$ on which $\rho$ enjoys \eqref{1dim-eq}.
Let
\[ T(x)=\exp_x \big( \nabla \phi(x) \big) =x+\phi'(x) \]
be the unique optimal transport map from
$\omega^{\ell}$ to $\mu^{\ell}$ by the $1$-dimensional version of Theorem~\ref{th:FG},
thus $\phi:\ell \longrightarrow \R$ is locally Lipschitz and twice differentiable almost everywhere.
For $t \in [0,1]$, set 
\[
T_t(x):=\exp_x \big( t\nabla \phi(x) \big) =x+t\phi'(x),
 \qquad \mu^{\ell}_t:=(T_t)_\sharp \omega^{\ell} =\rho^{\ell}_t \omega^{\ell}.
\]
Note that $W_2(\mu^{\ell}_s,\mu^{\ell}_t)=|s-t|W_2(\mu^{\ell}_0,\mu^{\ell}_1)$ for $s,t \in [0,1]$.

On the one hand,
we deduce from the Talagrand inequality \eqref{Talagrand}
and the $\CD(K,\infty)$-inequality \eqref{eq:CD} that, for any $t \in [0,1]$,
\begin{equation}\label{cvx}
\frac{K}{2}W_2(\mu^{\ell}_0,\mu^{\ell}_t)^2 \leq H_{\omega^{\ell}}(\mu^{\ell}_t)
 \leq (1-t) H_{\omega^{\ell}}(\mu^{\ell}_0)+t H_{\omega^{\ell}}(\mu^{\ell}_1)
 -\frac{K}{2}t(1-t)W_2(\mu^{\ell}_0,\mu^{\ell}_1)^2.
\end{equation}
On the other hand, $H_{\omega^{\ell}}(\mu^{\ell}_0)=H_{\omega^{\ell}}(\omega^{\ell})=0$
and the hypothesis \eqref{1dim-eq} yield that the RHS of \eqref{cvx} coincides with
\[
\frac{K}{2} t^2 W_2(\mu^{\ell}_0,\mu^{\ell}_1)^2 =\frac{K}{2} W_2(\mu^{\ell}_0,\mu^{\ell}_t)^2.
\]
Therefore all the inequalities in \eqref{cvx} become equality.

Equality in the $\CD(K,\infty)$-inequality implies rigidity-type results
(see \cite{BK} for a related work).
To this end, we have a closer look on the derivation of the $\CD(K,\infty)$-inequality \eqref{eq:CD}
from $\Psi''_{\ell} \ge K$ in Theorem~\ref{needle}(\ref{needle2})
(recall $d\omega^{\ell} =\e^{-\Psi_{\ell}} \,dx$).
Let us consider
\[ J_t(x) :=T_t'(x) \e^{\Psi_{\ell}(x)-\Psi_{\ell}(T_t(x))}
 = \big( 1+t\phi''(x) \big) \e^{\Psi_{\ell}(x)-\Psi_{\ell}(T_t(x))} \]
existing almost everywhere.
Recall from \eqref{MAonM} that
\begin{equation}\label{eq:MA}
\rho^{\ell}_t\big( T_t(x) \big) J_t(x) =\rho^{\ell}_0(x)
\end{equation}
holds almost everywhere.
We deduce from $\Psi''_{\ell} \ge K$ and $|\phi'(x)|=d_M(x,T(x))$ that
\begin{equation}\label{Jconcave}
\frac{d^2[\log J_t(x)]}{dt^2}
 =-\Psi''_{\ell} \big( T_t(x) \big) \phi'(x)^2 -\frac{\phi''(x)^2}{(1+t\phi''(x))^2}
 \le -Kd_M\big (x,T(x) \big)^2
\end{equation}
in the weak sense.
Since
\[ H_{\omega^{\ell}}(\mu^{\ell}_t) =\int_{\ell} \log \rho^{\ell}_t \,d\mu^{\ell}_t
 =\int_{\ell} \log \rho^{\ell}_t (T_t) \,d\omega^{\ell}
 =\int_{\ell} (\log \rho^{\ell}_0 -\log J_t) \,d\omega^{\ell} \]
by $\mu^{\ell}_t=(T_t)_{\sharp} \omega^{\ell}$ and \eqref{eq:MA},
we obtain from \eqref{Jconcave} the $\CD(K,\infty)$-inequality
\[ H_{\omega^{\ell}}(\mu^{\ell}_t) \le (1-t)H_{\omega^{\ell}}(\mu^{\ell}_0) +tH_{\omega^{\ell}}(\mu^{\ell}_1)
 -\frac{K}{2}t(1-t) \int_{\ell} d_M\big( x,T(x) \big)^2 \,\omega^{\ell}(dx). \]
Now we find from \eqref{Jconcave} that the equality \eqref{1dim-eq} shows $\phi''(x)=0$,
hence $\phi$ is an affine function on $\ell$,
as well as $\Psi''_{\ell}(T_t(x))d_M(x,T(x))^2 =Kd_M(x,T(x))^2$ for all $t$.

If $\phi' \neq 0$, then $T$ is a translation and hence $\ell$ is isometric to the real line $\R$.
Moreover, by $\Psi''_{\ell} \equiv K$,
$\omega^{\ell}$ is a Gaussian measure on $\ell$ with variance $K^{-1}$.
Notice also that, since the restriction of $\rho$ on $\ell$
is the ratio between two Gaussian densities with the same variance,
$\log\rho$ is also an affine function on $\ell$.

If $\phi' \equiv 0$, then $T$ is identity and $\mu^{\ell}=\omega^{\ell}$,
thereby $\rho \equiv 1$ on $\ell$.
In both cases $\rho>0$ holds, therefore we have $\supp\, \mu=M$.

\subsection{Integration in needles}

In the previous step we showed that $h:=\log\rho$ on each $\ell$
is either an affine function on the $1$-dimensional Gaussian space or constant.
In both cases $h$ yields equality of the Poincar\'e inequality \eqref{eq:Poin} on $\ell$,
namely
\[ \int_\ell h^2 \,d\omega^{\ell} -\bigg( \int_\ell h \,d\omega^{\ell} \bigg)^2
 =\frac{1}{K} \int_\ell |\nabla^{\ell} h|^2 \,d\omega^{\ell}. \]
Recall from \eqref{1dim-eq} that $|\nabla^\ell \rho|=|\nabla \rho|$ holds almost everywhere.
Therefore, by Theorem~\ref{needle}(\ref{needle1}) and the Cauchy--Schwarz inequality
on $(L,\nu)$,
\begin{align*}
\int_M h^2 \,d\omega -\bigg( \int_M h \,d\omega \bigg)^2
&= \int_L \int_\ell h^2 \,d\omega^{\ell} \,\nu(d\ell)
 -\bigg( \int_L \int_\ell h \,d\omega^{\ell} \,\nu(d\ell) \bigg)^2 \\
&\ge \int_L \int_\ell h^2 \,d\omega^{\ell} \,\nu(d\ell)
 -\int_L \bigg( \int_\ell h \,d\omega^{\ell} \bigg)^2 \,\nu(d\ell) \\
&= \frac{1}{K} \int_L \int_\ell |\nabla^{\ell} h|^2 \,d\omega^{\ell} \,\nu(d\ell) \\
&= \frac{1}{K} \int_M |\nabla h|^2 \,d\omega.
\end{align*}
Hence $h$ gives equality of the Poincar\'e inequality on $M$.
Thanks to \cite[Theorem~2]{CZ},
$(M,g,\omega)$ splits off the $1$-dimensional Gaussian space
as described in Theorem~\ref{main}(\ref{main1}).
Moreover, by the construction in~\cite{CZ},
every $\ell$ is Gaussian and $|\nabla h|$ is constant on $M$.
Since $\mu \neq \omega$ ensures $|\nabla h| \neq 0$,
$\mu$ is a translation of $\omega$ along the $\R$-direction.
This completes the proof of Theorem \ref{main}(\ref{main2}).

\section{Further problems}\label{outro}

Theorem~\ref{main} is concerned with one of the simplest situations
of the quantitative study of geometric or functional inequalities,
there are many related open problems.

\begin{enumerate}[(A)]
\item
It is natural to consider the analogue of Theorem~\ref{main} for $N \in [n,\infty)$.
In this case, however, there is no extremal function.
Precisely, for instance, tracing the proof of the dimensional logarithmic Sobolev inequality
under $\Ric_N \ge K>0$,
\[ H_{\omega}(\rho \omega) \le \frac{N-1}{2KN} I_{\omega}(\rho \omega), \]
in \cite[\S 5.7]{BGL}, we see that equality forces $\rho$ to be constant.
Therefore we can only assume the existence of a sequence of functions $\rho_i$ such that
\[ \lim_{i \to \infty} \frac{H_{\omega}(\rho_i \omega)}{I_{\omega}(\rho_i \omega)}
 =\frac{N-1}{2KN}. \]
According to the argument in \cite{OV} (see also \cite[\S 6]{LV}),
we may expect that $h_i:=\rho_i-1$ provides the sharp spectral gap, namely
\[ \lim_{i \to \infty} \frac{1}{\|h_i\|_{L^2(\omega)}^2} \int_M |\nabla h_i|^2 \,d\omega
 =\frac{KN}{N-1}. \]
To this end, we need to analyze the asymptotic behavior of $h_i$.
We remark that it is unclear if the needle decomposition is useful in this strategy.

\item
Cheng--Zhou's characterization in \cite{CZ} of the sharp spectral gap
was extended to $\RCD(K,\infty)$-spaces in \cite{GKKO}
($\RCD(K,\infty)$ stands for the \emph{Riemannian curvature-dimension condition}
corresponding to $\Ric_{\infty} \ge K$, for which we refer to \cite{AGS,EKS}).
Hence it is natural to expect the analogue to Theorem~\ref{main} for $\RCD(K,\infty)$-spaces.
Currently, however, the needle decomposition is known only for $\RCD(K,N)$-spaces with $N \in [1,\infty)$
(or, more generally, essentially non-branching $\CD(K,N)$-spaces), see \cite{CM1}.
In order to generalize Theorem~\ref{main},
we may generalize the needle decomposition to $\RCD(K,\infty)$-spaces,
or modify the discussion in this article to directly show that $\log\rho$
attains the sharp spectral gap (then \cite{GKKO} applies).

\item
It is also possible to consider a generalization of Theorem~\ref{main} for Finsler manifolds,
where the curvature-dimension condition and the needle decomposition are available 
(\cite{Oint,Oneedle}).
In this case, however, the splitting phenomenon is not as simple as the Riemannian world,
see \cite{Osplit} for the case of Cheeger--Gromoll splitting theorem.

\item
One of the most challenging problems is a quantitative version of
the logarithmic Sobolev inequality,
namely analysis of the situation where we have
\[ \int_M \rho \log \rho \,d\omega
 \ge \frac{1}{2K'} \int_M \frac{|\nabla \rho|^2}{\rho} \,d\omega \]
for some $K'>K$ close to $K$.
Quantitative properties of geometric or functional inequalities on Riemannian manifolds
are less investigated than the Euclidean or Gaussian cases.
The recent important work \cite{CMM} on a quantitative isoperimetric inequality
on $\CD(K,N)$-spaces with $K>0$ and $N \in (1,\infty)$ made use of the needle decomposition.
On Euclidean spaces equipped with Gaussian measures,
quantitative versions are studied for the isoperimetric inequality in \cite{CFMP,MN,El,BBJ}
and for the logarithmic Sobolev inequality in \cite{BGRS,FIL,CF}.

\item
Another interesting problem is the rigidity of gradient estimates,
in this case $K$ can be nonpositive as well.
See \cite{ABS} for a recent result in this direction.
\end{enumerate}
\medskip

{\it Acknowledgements}.
SO was supported in part by JSPS Grant-in-Aid for Scientific Research (KAKENHI) 15K04844, 17H02846.
AT was supported in part by JSPS Grant-in-Aid for Scientific Research (KAKENHI) 15K17536, 16KT0132.

{\small

}

\end{document}